\documentclass[10pt]{article}

\usepackage{amsmath}
\usepackage{mathtools}
\usepackage{upquote}
\usepackage{physics}
\usepackage{comment}

\usepackage[utf8]{inputenc}
\usepackage{arxiv}
\usepackage{graphicx}
\usepackage{multicol}
\usepackage{threeparttable}
\usepackage{subcaption}
\usepackage[none]{hyphenat}
\usepackage{titlesec}
\renewcommand{\thesection}{\Roman{section}}
\titleformat{\section}
  {\normalfont\normalsize\scshape}{\thesection.}{1em}{}
\renewcommand{\thesubsection}{\Alph{subsection}}
\titleformat{\subsection}
  {\normalfont\normalsize\itshape}{\thesubsection.}{1em}{}
\titleformat{\subsubsection}
  {\normalfont\normalsize\itshape}{\thesubsection.\thesubsubsection.}{1em}{}
\renewcommand{\thesubsubsection}{\roman{subsubsection}}
\usepackage{titling}
\predate{}
\postdate{}

\graphicspath{ {./Images/} }

\author{
\textbf{Jayanth S. Pratap} \\
  Dulles Math and Science Academy\\\\
  Presented at the 3\textsuperscript{rd} Systems Approaches to Cancer Biology Conference\\November 11\textsuperscript{th} - 13\textsuperscript{th}, 2020.
}
\title{\textsc{An Optimal Control Strategy for Mathematically Modeling Cancer Combination Therapy}}
\date{}

\begin{document}

\maketitle

\begin{abstract}

    While the use of combination therapy is increasing in prevalence for cancer treatment, it is often difficult to predict the exact interactions between different treatment forms, and their synergistic/antagonistic effects on patient health and therapy outcome. In this research, a system of ordinary differential equations is constructed to model nonlinear dynamics between tumor cells, immune cells, and three forms of therapy: chemotherapy, immunotherapy, and radiotherapy. This model is then used to generate optimized combination therapy plans using optimal control theory. In-silico experiments are conducted to simulate the response of the patient model to various treatment plans. This is the first mathematical model in current literature to introduce radiotherapy as an option alongside immuno- and chemotherapy, permitting more flexible and effective treatment plans that reflect modern therapeutic approaches.

\end{abstract}

\begin{multicols}{2}

\section{Introduction}

Combination therapy in cancer treatment is currently a highly active area of research. Immunotherapy, chemotherapy, and radiotherapy are often used together in some combination to maximize the likelihood of tumor regression, but in these cases it is of paramount importance to consider both the synergistic and antagonistic effects of these drugs when combined. Chemotherapy and radiotherapy often have toxic effects on the patient's immune system, and immunotherapy can be used to counteract these effects by stimulating immune system activation.

Quantitative modeling provides a promising avenue for analyzing and optimizing the synergistic effects of different therapies. In this project, the first phase consists of constructing a model of ordinary differential equations to capture interactions between tumor cells, immune effector cells, and three forms of therapeutic intervention---immuno-, chemo-, and radiotherapy. Biological parameters fitted to murine trials and patient data are used to simulate the model. The second phase is setting up different optimal control scenarios in fixed final-time form to minimize different quantities over the course of treatment, such as average tumor burden, final tumor burden, and therapy delivered, as well as introducing a term to maximize patient health in the form of average immune cell population.

\section{ODE Model Formulation}

The model is comprised of a two-dimensional system, and is modeled around the Lotka-Volterra predator-prey system, with form
\begin{equation*}
    \dot{x} = \alpha x - \beta xy
\end{equation*}
\begin{equation*}
    \dot{y} = \delta xy - \gamma y,
\end{equation*} 

where $x$ represents the prey species and $y$ represents the predator species.

In the system of tumor-immune dynamics, the predator species $y$ is considered as the effector cell population (immune cells with the ability to kill tumor cells), and the prey species $x$ is represented as the tumor cell population.

\subsection{ODE Model: General Form of Equations}

There are five populations whose dynamics are considered in the ODE model:
\begin{itemize}
    \item $N(t)$, tumor cell population
    \item $E(t)$, effector cell population
    \item $I(t)$, immunotherapy dose level
    \item $C(t)$, chemotherapy dose level
    \item $R(t)$, radiotherapy dose level
\end{itemize}

Furthermore, there are several forms of terms that must be included in the equations:
\begin{itemize}
    \item Net growth term ($G_N$, $G_E$, $G_I$, $G_C$, $G_R$)
    \item Tumor-immune agonism/antagonism ($A_{E\rightarrow N}$, $A_{N\rightarrow E}$)
    \item Chemo-therapeutic effects ($C_N$, $C_E$)
    \item Radio-therapeutic effects ($R_N$, $R_E$)
    \item Immuno-therapeutic effects ($I_E$)
    \item Therapy intervention ($u_I$, $u_C$, $u_R$)
\end{itemize}

\subsubsection{Growth and death terms}

In the absence of therapy, tumor cells ($N$) grow in some fashion limited by \textit{in vivo} carrying capacity. The limited-growth model used in describing this is logistic growth, based on murine data. This gives us the form
\begin{equation*}
    G_{N} = \alpha N(1 - \beta N).
\end{equation*}

Immune cells are assumed to be generated at a constant rate, and to have a natural lifespan until death. This gives us the form 
\begin{equation*}
    G_{E} = \zeta - \lambda E.
\end{equation*}

For all three forms of therapy analyzed, we assume an exponential decay function adequately models the elimination of the dose from the body, thus
\begin{equation*}
    G_{I} = -\tau E,
\end{equation*}
\begin{equation*}
    G_{C} = -\omega C, \text{ and }
\end{equation*}
\begin{equation*}
    G_{R} = -\chi R.
\end{equation*}

\subsubsection{Tumor-immune agonism/antagonism}

Tumor cells are assumed to be deactivated by circulating effector cells in a manner that scales with the probability of an effector/tumor cell encounter, thus giving the form
\begin{equation*}
    A_{E\rightarrow N} = - \gamma NE.
\end{equation*}

There is an immunostimulatory effect as effector cells come into contact with tumor cells. This is due to cytokine-induced activation, where chemokines (immune signaling molecules) and inflammatory markers induce the increased production and activation of effector cells.
This term is modeled with Michaelis-Menten kinetics because immunostimulation asymptotically approaches a maximum rate due to limitation in the rate of immune cell circulation and cell signaling. Thus we use form
\begin{equation*}
    A_{N\rightarrow E} = \frac{\eta EN}{\theta + N}
\end{equation*} where $\eta$ represents maximum rate of activation and $\theta$ represents the half-saturation constant, equivalent to the value of tumor cells for which the activation rate is half of its maximum.

\subsubsection{Therapeutic effects}

The form of immunotherapy analyzed in this project is IL-2 (interleukin-2) immunotherapy; IL-2 is a naturally occurring cytokine in the body, and its effect on immune system activation is thus described in the same form as $M_{N\rightarrow E}$:
\begin{equation*}
    I_{E} = \frac{\nu EI}{\rho + I}.
\end{equation*}

Chemotherapy kill terms reflect fractional cell kill: due to chemotherapeutic agents such as doxorubicin, cisplatin, and bleomycin being mostly or exclusively effective towards tumor cells undergoing specific phases of the cell cycle, there is always only a particular fraction of cells that are killed with each dosage, proportional to the dosage itself. Additionally, the dose-response curves of chemotherapeutic drugs often demonstrate bounded efficacy of chemotherapy, and thus a saturation term is used to give near-linear kill rate for low drug concentrations, with plateauing kill rate for higher concentrations. With these considerations, we use terms
\begin{equation*}
    I_{N} = - \delta N(1-{\rm e}^{-C}) \text{ and }
\end{equation*}
\begin{equation*}
    C_{E} = - \mu E(1-{\rm e}^{-C}).
\end{equation*}

Radiotherapy kill terms reflect the linear-quadratic model, used in radiation biology to reflect the phenomenon that cell survival fractions after irradiation often take the form of an exponential function with a linear and quadratic term. We use terms
\begin{equation*}
    R_{N} = - N(\epsilon R + \kappa R^2) \text{ and }
\end{equation*}
\begin{equation*}
    R_{E} = - E(\sigma R + \phi R^2).
\end{equation*}

The ratios $\frac{\epsilon}{\kappa}$ and $\frac{\sigma}{\phi}$, as per clinical trial data in \cite{mehta_2004}, are set as $\frac{\epsilon}{\kappa} = \frac{\sigma}{\phi} = 10$. 

\subsubsection{Therapy intervention}

We allow functions $u_{I}(t)$, $u_{C}(t)$, $u_{R}(t)$ to model the infusion rate of immuno-, chemo-, and radio-therapy respectively, defined in terms of time. These are intended to serve as controls to optimize treatment, which will be examined later in the project.

\subsection{ODE Model: Specific Forms of Equations}

Considering the developments in the previous section, the ODE system has been written as follows: 

\end{multicols}

\begin{subequations}
\begin{equation}
    \dot{N}(t) =
    \alpha N(1 - \beta N)
    - \gamma NE
    - \delta N(1-e^{-C})
    - N(\epsilon\hat{R} + \kappa\hat{R}^2)
\end{equation}
\begin{equation}
    \dot{E}(t) =
    \zeta - \lambda E
    + \frac{\eta EN}{\theta + N}
    + \frac{\nu EI}{\rho + I}
    - \mu E(1-e^{-C})
    - E(\sigma R + \phi R^2)
\end{equation}
\begin{equation}
    \dot{I}(t) = u_{I}(t) - \tau I
\end{equation}
\begin{equation}
    \dot{C}(t) = u_{C}(t) - \omega C
\end{equation}
\begin{equation}
    \dot{R}(t) = u_{R}(t) - \chi R
\end{equation}
\end{subequations}

\subsection{Model Parameters}

The parameters used in the project are sourced from various past models which have analyzed data from murine experiments and human clinical trials for curve-fitting. See below for a table of parameters, their descriptions, and the values used.

\begin{center}
\begin{scriptsize}
\begin{threeparttable}
    \caption{Estimated Model Parameter Values}
    \begin{tabular}{l*{4}{c}r}
        \\ \hline \\
        Parameter & Description & Units & Estimated Value & Source\\\\
        \hline \\
         $\alpha$ & Maximum growth rate of tumor cells & day$^{-1}$ & 1.80 $\times$ 10$^{-1}$ & \cite{kuznetsov}\\
         $\beta$ & Reciprocal of carrying capacity of tumor cells & cells$^{-1}$ & 2.00 $\times$ 10$^{-9}$ & \cite{kuznetsov}\\
         $\gamma$ & Effector-cell-induced tumor death rate & cells$^{-1}$ day$^{-1}$ & 1.101 $\times$ 10$^{-7}$ & \cite{kuznetsov}\\
         $\delta$ & Chemotherapy (Bleomycin) kill rate coefficient for tumor cells & day$^{-1}$ & 9.00 $\times$ 10$^{-1}$ & \cite{pillis_gu_radunskaya_2006}\\
         $\epsilon$ & Radiotherapy linear kill coefficient for tumor cells & Gy$^{-1}$ & 3.98 $\times$ 10$^{-2}$ & \cite{predictioncombined}\\
         $\kappa$ & Radiotherapy quadratic kill coefficient for tumor cells & Gy$^{-2}$ & 3.98 $\times$ 10$^{-3}$ & \cite{predictioncombined}\footnotemark\\
         $\zeta$ & Constant effector cell production rate & cells day$^{-1}$ & 1.30 $\times$ 10$^{4}$ & \cite{kuznetsov}\\
         $\eta$ & Maximum cytokine activation rate of effector cells & day$^{-1}$ & 1.245 $\times$ 10$^{-1}$ & \cite{kuznetsov}\\
         $\theta$ & Half-saturation constant for cytokine activation of effector cells & cells & 2.019 $\times$ 10$^{7}$ & \cite{kuznetsov}\\
         $\lambda$ & Death rate of effector cells & day$^{-1}$ & 4.12 $\times$ 10$^{-2}$ & \cite{kuznetsov}\\
         $\nu$ & Maximum immunotherapy (IL-2) activation rate of effector cells & day$^{-1}$ & 1.245 $\times$ 10$^{-1}$ & \cite{kirschner_panetta_1998}\\
         $\rho$ & Half-saturation constant for immunotherapy (IL-2) activation of effector cells & cells & 2.00 $\times$ 10$^{7}$ & \cite{kirschner_panetta_1998}\\
         $\mu$ & Chemotherapy (Bleomycin) kill rate coefficient for effector cells & day$^{-1}$ & 6.00 $\times$ 10$^{-1}$ & \cite{pillis_gu_radunskaya_2006}\\
         $\sigma$ & Radiotherapy linear kill coefficient for effector cells & Gy$^{-1}$ & 3.98 $\times$ 10$^{-2}$ & Assumed\footnotemark\\
         $\phi$ & Radiotherapy quadratic kill coefficient for effector cells & Gy$^{-2}$ & 3.98 $\times$ 10$^{-3}$ & Assumed\footnotemark\\
         $\tau$ & Decay rate of immunotherapy (IL-2) & day$^{-1}$ & 1.00 $\times$ 10$^{1}$ & \cite{kirschner_panetta_1998}\\
         $\omega$ & Decay rate of chemotherapy (Bleomycin) & day$^{-1}$ & 9.00 $\times$ 10$^{-1}$ & \cite{pillis_gu_radunskaya_2006}\\
         $\chi$ & Decay rate of radiotherapy & day$^{-1}$ & 1.1 $\times$ 10$^{-2}$ & \cite{mkango_shaban_mureithi_ngoma_2019}\\
         $I_{max}$ & Maximum tolerable dose of immunotherapy (IL-2) & IU & 7.20 $\times$ 10$^{5}$ & \cite{dutcher_schwartzentruber_kaufman_agarwala_tarhini_lowder_atkins_2014}\\
         $C_{max}$ & Maximum tolerable dose of chemotherapy (Bleomycin) & IU & 3.00 $\times$ 10$^{4}$ & \cite{aston_hope_nowak_robinson_lake_lesterhuis_2017}\\
         $R_{max}$ & Maximum tolerable dose of radiotherapy & Gy & 4.50 $\times$ 10$^{1}$ & \cite{mehta_2004}\\\\
         \hline \\
    \end{tabular}
    \begin{tablenotes}
        \item [1, 3] Assumed based on ratio from \cite{mehta_2004}.
        \item [2] Assumed as equal to tumor cell kill coefficient as per \cite{mkango_shaban_mureithi_ngoma_2019}. 
        \\
    \end{tablenotes}
\end{threeparttable}\bigskip
\end{scriptsize}
\end{center}

\begin{multicols}{2}

\subsubsection{Non-dimensionalization of Parameters}

We define the nondimensionalized state variables as follows: \\
\begin{align*}
    \hat{N} &= N/N_0, \quad
    &\hat{E} = E/E_0, \\
    \hat{I} &= I/I_0, \quad
    &\hat{C} = C/C_0, \\
    \hat{R} &= R/R_0,
\end{align*} \\ where the scaling factors are chosen such that \\
\begin{align*}
    N_0 &= 10^6, \quad
    &E_0 = 10^6, \\
    I_0 &= I_{max}, \quad
    &C_0 = C_{max}, \\
    R_0 &= R_{max},
\end{align*} \\to rescale the populations by a reasonable order of magnitude.

Time is rescaled with respect to tumor cell deactivation such that 
\begin{equation*}
    t' = t/t_0, \quad \text{where}\quad    t_0 = (\gamma N_0)^{-1},
\end{equation*} the time scaling factor. This eliminates the parameter $\gamma$.

The parameter values are non-dimensionalized with the following transformations:
\small
\begin{align*}
    &\alpha^\ast = \alpha \cdot t_0, \quad
    \beta^\ast = \beta \cdot N_0, \quad
    \gamma^\ast = \gamma \cdot N_0 \cdot t_0, \quad
    \delta^\ast = \delta \cdot t_0, \\
    &\epsilon^\ast = \epsilon \cdot t_0, \quad
    \kappa^\ast = \kappa \cdot t_0, \quad
    \zeta^\ast = \frac{\zeta \cdot t_0}{E_0}, \quad
    \lambda^\ast = \lambda \cdot t_0, \\
    &\eta^\ast = \frac{\eta}{E_0}, \quad
    \theta^\ast = \theta \cdot t_0, \quad
    \nu^\ast = \nu \cdot t_0, \quad
    \rho^\ast = \frac{\rho}{I_0}, \\
    &\mu^\ast = \mu \cdot t_0, \quad
    \sigma^\ast = \sigma \cdot R_0, \quad
    \phi^\ast = \phi \cdot R_{0}^2, \\
    &\tau^\ast = \tau \cdot t_0, \quad
    \omega^\ast = \omega \cdot t_0, \quad
    \chi^\ast = \chi \cdot t_0.
\end{align*} \normalsize Dropping the stars for notational clarity, the final non-dimensionalized system is now given by:

\end{multicols}

\begin{subequations}
\begin{equation}
    \dot{\hat{N}}(t') =
    \alpha\hat{N}(1 - \beta\hat{N})
    - \hat{N}\hat{E}
    - \delta\hat{N}(1-e^{-\hat{C}})
    - \hat{N}(\epsilon\hat{R} + \kappa\hat{R}^2)
\end{equation}
\begin{equation}
    \dot{\hat{E}}(t') =
    \zeta - \lambda\hat{E}
    + \frac{\eta\hat{E}\hat{N}}{\theta + \hat{N}}
    + \frac{\nu\hat{E}\hat{I}}{\rho + \hat{I}}
    - \mu\hat{E}(1-e^{-\hat{C}})
    - \hat{E}(\sigma\hat{R} + \phi\hat{R}^2)
\end{equation}
\begin{equation}
    \dot{\hat{I}}(t') = u_{\hat{I}}(t') - \tau\hat{I}
\end{equation}
\begin{equation}
    \dot{\hat{C}}(t') = u_{\hat{C}}(t') - \omega\hat{C}
\end{equation}
\begin{equation}
    \dot{\hat{R}}(t') = u_{\hat{R}}(t') - \chi\hat{R}
\end{equation}
\end{subequations}

\vspace{5mm}

\begin{multicols}{2}

\subsection{Treatment-Free Equilibria and Stability}

In this section we consider the equilibrium behavior of the system in the absence of treatment. For convenience of notation, let $x_1 \equiv \hat{N}$ and $x_2 \equiv \hat{E}$. First we must find the nullclines of this system, or the curves along which $\dv{x_1}{t} = \dv{x_2}{t} = 0$. The intersection of these nullclines gives the fixed points of the system. 

Setting $\dv{x_1}{t} = 0$ gives two nullclines of the system, one of which represents a tumor-free equilibrium and the other a non-zero tumor population at equilibrium:

\begin{subequations}
    \begin{equation}
        x_1 = 0
    \end{equation}
    \begin{equation}
        x_2 = \alpha(1 - \beta x_1) \equiv h(x_1)
    \end{equation}
\end{subequations}

Setting $\dv{x_2}{t} = 0$ gives the last nullcline:

\begin{equation}
    x_2 = \frac{\zeta}{\lambda - \frac{\eta x_1}{\theta + x_1}} \equiv j(x_1)
\end{equation}

Using these nullclines, we can now analyze the various equilibrium states of the system.

\subsubsection{Tumor-Free Equilibrium}

The intersection of $h(x_1)$ with $x_1 = 0$ gives a fixed point with coordinates $(x_1, x_2) = (0, \frac{\zeta}{\lambda}) \equiv \xi_1$.

To permit stability analysis, we linearize about point $\xi_1$. Let $f_1(x_1, x_2)$ represent $\dv{x_1}{t}$, and $f_2(x_1, x_2)$ represent $\dv{x_2}{t}$ for ease of notation. The Jacobian matrix for this system is of the form

\begin{equation*}
\mathbf{J}(x_1, x_2) =
\begin{bmatrix}
  \frac{\partial f_1}{\partial x_1} & 
    \frac{\partial f_1}{\partial x_2} \\[1ex]
  \frac{\partial f_2}{\partial x_1} & 
    \frac{\partial f_2}{\partial x_2} \\[1ex]
\end{bmatrix}
\end{equation*}

Evaluating at $\xi_1$ gives a Jacobian of

\begin{equation*}
\mathbf{J}(\xi_1) =
\begin{bmatrix}
   \alpha - \frac{\zeta}{\lambda} & 
    0 \\[1ex]
  \frac{\eta\zeta}{\theta\lambda}& 
    - \lambda\\[1ex]
\end{bmatrix}
\end{equation*}

By the Hartman-Grobman theorem, the real components of both eigenvalues of the Jacobian matrix must be negative for the given equilibrium state to be stable. Eigenvalues are

\[ \lambda_{1} = \alpha - \frac{\zeta}{\lambda}\]
\[ \lambda_{2} = -\lambda < 0\]

Thus, for the tumor-free fixed point $\xi_1$ to be a stable equilibrium point for the system, it holds that

\begin{equation*}
\alpha < \frac{\zeta}{\lambda}
\end{equation*}

This allows us to see that to eliminate the tumor cell population at the stable equilibrium of the system, the maximum rate of tumor proliferation must be less than the equilibrium tumor cell population. If this condition is not met, the tumor-free equilibrium will be unstable.

\subsubsection{Nonzero Equilibria}

There may be zero, one, or two additional equilibria besides $\xi_1$, depending on the relationship between $h(x_1)$ and $j(x_1)$.

Intersecting $h(x_1)$ and $j(x_1)$ gives a quadratic of the following form, whose solutions are values of $x_1$ for all other equilibria of the system.

\begin{equation}
    Ax_1^2 + Bx_1 + C = 0,
\end{equation}

where

\begin{align*}
    & A = \beta(\lambda - \eta), \quad B = \frac{\zeta}{\alpha} + \beta\lambda\eta + \eta - \lambda, \\
    & \text{and} \quad C = \theta\Big(\frac{\zeta}{\alpha} - \lambda\Big).
\end{align*}

The number of solutions to this quadratic depends on values of $A$, $B$, and $C$, thus depending on the parameter values being used. For the chosen parameter values, the quadratic formula yields two positive values for $x_1$ at equilibrium. However, one of these values yields a negative $x_2$ and is thus eliminated. Therefore, in this case our analysis will only include one additional nonzero fixed point $\xi_2$.

\subsection{Limit Cycles}

The Dulac-Bendixson criterion can be used to prove the absence of closed orbits in this system. The Dulac-Bendixson criterion states that if there are no closed orbits, there exists a function $\phi$ such that

\begin{align}
    \frac{\partial}{\partial x_1} (\phi \dot{x_1}) + \frac{\partial}{\partial x_2} (\phi \dot{x_2}) \neq 0
\end{align}

We choose $\phi = \frac{1}{x_1 x_2}$. 

\begin{align}
    &\frac{\partial}{\partial x_1} (\phi \dot{x_1}) + \frac{\partial}{\partial x_2} (\phi \dot{x_2}) \\
    & \quad = -(\frac{\alpha \beta}{x_2} + \frac{\zeta}{x_1 x_{2}^2}) \neq 0.
\end{align}

Therefore, the Dulac-Bendixson criterion is satisfied and it holds that there are no limit cycles or homoclinic connections in the system. Evaluating the criterion for both treatment-free and treatment-present cases gives this conclusion. This proves that the system will always collapse to an equilibrium point, with or without any therapeutic intervention.

\subsection{Numerical Solution}

We numerically solve the ODE system using the Python Pyomo library. As mentioned before, the estimated parameter values are such that there are two steady states of the system: $\xi_1$ being an unstable tumor-free state and $\xi_2$ being a stable steady state with a non-zero, dormant tumor population. Graphs of dynamics are shown below with various initial conditions.

\end{multicols}

\begin{figure}[h]
    \centering
    \includegraphics[width=16cm]{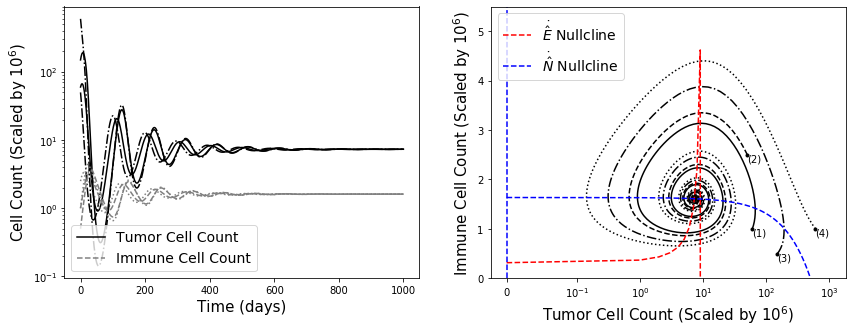}
    \caption{\textit{Time dynamics (left) and phase portrait (right) of tumor cells and immune cells without treatment. Initial conditions $(N, E) = (6 \times 10^{7}, 1 \times 10^{6}), (5 \times 10^{7}, 2.5 \times 10^{6}), (1.5 \times 10^{8}, 1 \times 10^{6}), (6 \times 10^{8}, 1 \times 10^{6})$}}
\end{figure}

\begin{multicols}{2}

\section{Optimal Control}

A general optimal control problem is stated as the following: a system of ODEs is given that models the dynamics of a system in terms of state variables and control inputs, such as
\begin{align}
    \dot{\textbf{x}} = f(\textbf{x}, \textbf{u}, t),
\end{align}
where $f$ is a vector-valued function and $\textbf{x}, \textbf{u}$ are vectors that describe the state of the system and values of control inputs, respectively. The goal is to drive the system to a final state given initial states of the system, while minimizing an objective function $\mathcal{J}$ through the time duration $[t_0, t_f]$. There may also be additional constraints added that must be satisfied for this duration, called path constraints. In the following sections, we aim to define and solve an optimal control problem to optimize treatment.

\subsection{Defining Objective Function}

In this section, we explore a variant of the problem in which the final time $t_f$ is fixed at a particular value, known as a fixed-final time problem in optimal control theory.

The objective function is the function that must be minimized through the course of treatment. In this problem, there are several possible metrics that can be used to represent the goals for treatment duration and outcome.

\subsubsection{Incorporating Physician Preference}

There must be some margin provided in the model to allow for physician preference to be incorporated into optimization. We introduce weighting coefficients $w_1$ and $w_2$ to allow weighted preference for two different parts of the objective function. Thus, the objective function is in the form

\begin{equation}
    \mathcal{J} = w_1 \mathcal{K} + w_2 \mathcal{L},
\end{equation}

where $\mathcal{K}$ and $\mathcal{L}$ represent chosen expressions involving different quantities in the model that are to be minimized. The rest of this section will explore various choices for these expressions, and the impact they have on optimization results. Since there are two terms in the objective function, the results will be Pareto-optimal. 

\subsubsection{Maximizing Average Immune Cell Count}

By selecting weighting coefficients such that $w_2 < 0$, one can choose to maximize a quantity throughout treatment duration. Ideally, patient health should be maximized throughout treatment, and a good representative marker for patient health is immune cell count. Thus, we select $\mathcal{L}$ to be average immune cell count, in the form

\begin{equation}
    \mathcal{L} = \frac{1}{t_f} \int_{0}^{t_f} \hat{E}(t) dt.
\end{equation}

This is the $\mathcal{L}$ used for all optimization scenarios explored in the rest of the project. Note that with varying values of $w_2 < 0$, a physician can select an optimization setup that reflects the preferred approach for treatment, choosing a more cautious route with higher values of $|w_2|$, and a more aggressive treatment plan with lower values of $|w_2|$.

\subsubsection{Minimizing Final Tumor Cell Count}

One possible aim is to drive the tumor population to as low a value as possible at the final time of treatment. In this case, the value of $\mathcal{K}$ would be defined as the tumor cell count at the end of treatment.

\begin{equation}
    \mathcal{K} = \hat{N}(t_f).
\end{equation}

\subsubsection{Minimizing Average Tumor Cell Count}

Another possible aim is to keep the tumor population as low as possible for not just the end, but the entire duration of treatment, to prevent development of malignancy. As such, $\mathcal{K}$ would be defined as the average tumor population over the treatment time.

\begin{equation}
    \mathcal{K} = \frac{1}{t_f} \int_{0}^{t_f} \hat{N}(t) dt.
\end{equation}

\subsubsection{Minimizing Treatment}

There is a limited set of initial conditions and parameter values where a solution is possible such that $N(t_f) = 0$. In these cases, it is possible to explore the solution space and optimize for the solution wherein the least treatment is used. Minimizing the provided treatment serves to also minimize toxicity to the patient and cost of the therapy. In this case, $\mathcal{K}$ would be defined as the integral of all treatment options used over the duration of therapy:

\begin{equation}
    \mathcal{K} = \int_{0}^{t_f} u_I(t) + u_C(t) + u_R(t) dt.
\end{equation}

When using this objective function, there must also be a constraint placed on the optimization setup, such that

\begin{equation}
    \hat{N}(t_f) = 0.
\end{equation}

This restricts the use of this objective function to only the cases where it is possible to eliminate the tumor -- otherwise, the constraint is infeasible and one of the previous two objective functions must be used instead.

\subsubsection{Pontryagin's Maximum Principle}

Recall that $u_I, u_C,$ and $u_R$ represent the infusion rate of immuno-, chemo-, and radiotherapy, respectively, at a given time.

For future notation, we define $\textbf{u}$ to be the vector
$\begin{bmatrix}
    u_I & u_c & u_R
\end{bmatrix}^T$,
and $\textbf{x}$ to be
$\begin{bmatrix}
    x_1 & x_2 & x_3 & x_4 & x_5 
\end{bmatrix}^T$,
where $x_3, x_4, x_5$ represent $\hat{I}, \hat{C}, \hat{R}$ respectively, and as shown earlier, $x_1, x_2$ represent $\hat{N}, \hat{E}$ respectively. This allows us to capture the dynamics in Equations 2a-2e as $\dot{\textbf{x}} = f(\textbf{x}, \textbf{u}, t)$.

Now that we have the correct form for the dynamics, we can use Pontryagin's Maximum Principle. Pontryagin's Maximum Principle defines the first-order conditions necessary for a particular solution to the optimal control problem to be optimal.

The Hamiltonian of the system is defined as 

\begin{align}
    \mathcal{H} = \mathcal{J}(\textbf{u}, t) + \lambda^{T}(t) f(\textbf{x}, \textbf{u}, t).
\end{align}

$\lambda(t)$ is the vector of co-state (adjoint) variables, which vary with time. The Maximum Principle states that to find an optimal control, the Hamiltonian in this case must be minimized. The necessary conditions for optimality are as follows:

\begin{enumerate}
    \item $\dot{\textbf{x}} = \frac{\partial H}{\partial \lambda} = f(\textbf{x}, \textbf{u}, t)$ (State Equation)
    \item $\dot{\lambda} = - \frac{\partial H}{\partial \textbf{x}}$ (Costate Equation)
    \item $\frac{\partial H}{\partial \textbf{u}} = 0$ (Optimal Control)
    \item $0 \leq \textbf{u} \leq 1$ (Constraint on Control Inputs)
    \item $\textbf{x}(t_0) = \textbf{x}_0$ (Initial Conditions)
\end{enumerate}

\subsubsection{Numerical Solutions}

The Python Pyomo library \cite{hart2011pyomo} provides a robust framework for numerical optimization. The co-state variables are calculated to minimize the Hamiltonian and meet the aforementioned constraints. Below are shown the results after optimizing with the chosen constraints and objective function.

\end{multicols}

\begin{figure}[h]
    \centering
    \begin{subfigure}{\linewidth}
        \centering
        \includegraphics[width=12cm]{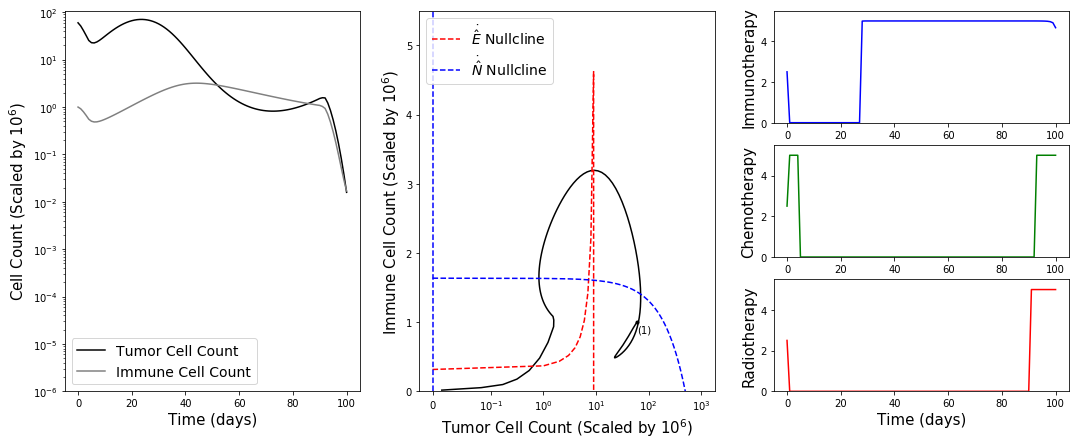}
        \caption{\textit{Minimize final tumor population.}}
    \end{subfigure}\par\medskip
    
    \begin{subfigure}{\linewidth}
        \centering
        \includegraphics[width=12cm]{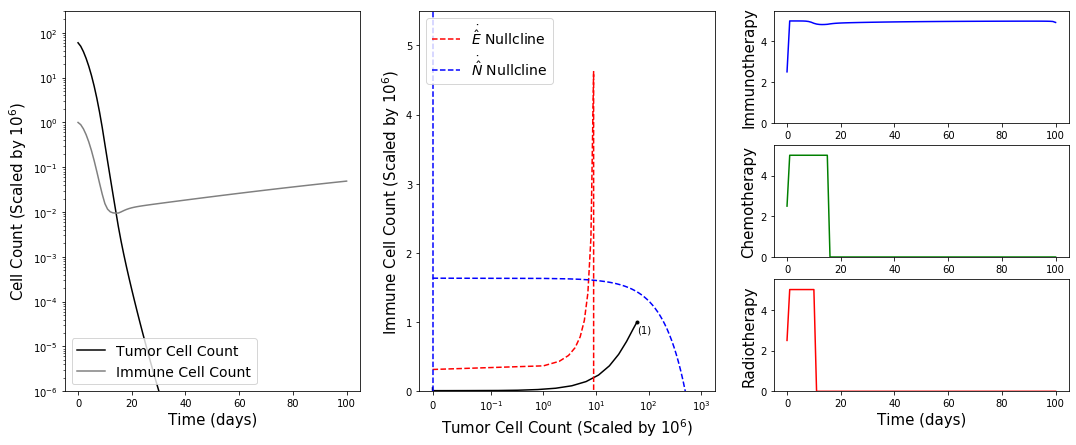}
        \caption{\textit{Minimize average tumor population.}}
    \end{subfigure}\par\medskip
    
    \begin{subfigure}{\linewidth}
        \centering
        \includegraphics[width=12cm]{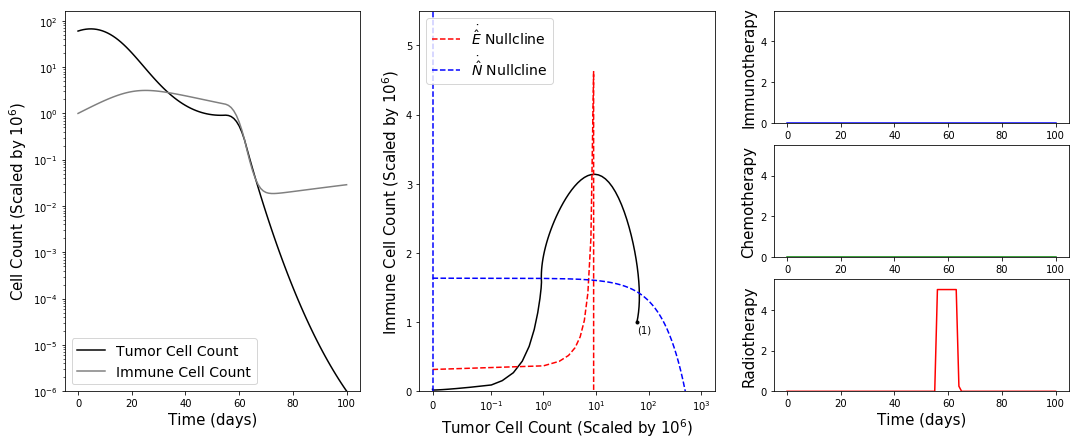}
        \caption{\textit{Minimize treatment administered, paired with constraint that final tumor population is zero.}}
    \end{subfigure}\par\medskip
    
    \caption{\textit{Time dynamics (left) and phase portrait (middle) of tumor cells and immune cells with treatment. Treatment regimen (right) showing the calculated control inputs for the duration of treatment. Initial conditions $(N, E) = (6 \times 10^{8}, 1 \times 10^{6})$. Weighting coefficients $w_1 = 1, w_2 = 0.25$.}}
\end{figure}

\begin{multicols}{2}

\subsection{Results and Discussion}

As shown in the computational experiments above, the optimal treatment regimen differs depending on the chosen objective function. When minimizing final tumor population, the optimal treatment plan involves peaks of bleomycin and radiotherapy near the last 10-day period, paired with high-dose IL-2 after the first month of treatment. When minimizing average tumor population, the regimen involves aggressive early peaks of bleomycin and radiotherapy for the first 2-week period and high-dose IL-2 for the entire duration. Finally, since eliminating the tumor is feasible for these parameters, the most optimal treatment regimen -- to eliminate the tumor with minimum impact to patient health -- is one where a single 10-day treatment cycle of radiotherapy is used. The dynamics will vary for different parameter values, due to the presence of additional equilibria states. Regardless, the model with the given parameters is an effective case study of generating an optimal treatment plan from pre-defined patient data values.

\section{Conclusion}

This project builds upon past contributions to computational oncology by extending the use of optimal control theory to combination therapy, as opposed to single-therapy treatment plans. This is a stepping stone to the possibility of personalized, optimized combination therapy regimens built around patient data. Future work will include optimization of treatment time and incorporating more biological phenomena into the model.

\nocite{villasana_radunskaya_2003}
\nocite{gluzman_scott_vladimirsky_2020}
\nocite{pillis_radunskaya_2003}

\bibliographystyle{acm}
\bibliography{bibliography}

\end{multicols}

\end{document}